\input amstex
\documentstyle{amsppt}
\linespacing{1.67}

\input xy
\xyoption{all}
\input graphicx.tex

\magnification=1200
\hsize=150truemm
\vsize=224.4truemm
\hoffset=4.8truemm
\voffset=12truemm

   \TagsOnRight
\NoBlackBoxes
\NoRunningHeads

\def\Square{\rlap{$\sqcup$}$\sqcap$}
\def\cqfd {\quad \hglue 7pt\par\vskip-\baselineskip\vskip-\parskip
{\rightline{\Square}}}

\redefine\O{{\Cal O}}
\define\T{{\Cal T}}
\define\M{{\Cal M}}
\define\C{{\Cal C}}
\define\A{{\Cal A}}
\define\R{{\Bbb R}}
\define\Z{{\Bbb Z}}
\define\Rt {${\R}$-tree}
\define\mi{^{-1}}
\let\thm\proclaim
\let\fthm\endproclaim
\let\inc\subset 
\let\ds\displaystyle
\let\ev\emptyset
\let\wtilde\widetilde
\let\ov\overline
\define\Aut{\text{\rm Aut}}
\define\MCG{\text{\rm PMCG}}
\define\Out{\text{\rm Out}}
\define\m{^{-1}}

\redefine\phi{\varphi}

\define\bbR{{\R}}
\define\ra{\rightarrow}
\def\emph#1{{\it #1\/}}
\define\rond#1{\overset{\circ} \to {#1}}%
\def\Tmin{T_{\text{min}}}

\newcount\tagno
\newcount\secno
\newcount\subsecno
\newcount\stno
\global\subsecno=1
\global\tagno=0
\define\ntag{\global\advance\tagno by 1\tag{\the\tagno}}

\define\sta{\the\secno.\the\stno
\global\advance\stno by 1}

\define\sect{\global\advance\secno by
1\global\subsecno=1\global\stno=1\
\the\secno. }

\define\subsect{\the\secno.\the\subsecno. \global\advance
\subsecno by 1}

\def\nom#1{\edef#1{\the\secno.\the\stno}}
\def\eqnom#1{\edef#1{(\the\tagno)}}

\newcount\refno
\global\refno=0
\def\nextref#1{\global\advance\refno by 1\xdef#1{\the\refno}}
\def\bref {\ref\global\advance\refno by 1\key{\the\refno}}

\nextref\AB
\nextref\BJ
\nextref\BF
\nextref\Bo
\nextref\BKM
\nextref\Ch
\nextref\Cl
\nextref\CuMo
\nextref\CV
\nextref\Fo
\nextref\FR
\nextref\GLP
\nextref\Gr
\nextref\Gu
\nextref\GLdef
\nextref\KW
\nextref\kmI
\nextref\KM
\nextref\Kr
\nextref\Le
\nextref\LL
\nextref\McC
\nextref\McM
\nextref\Pa
\nextref\Pab
\nextref\Remes
\nextref\Se
\nextref\SeIHES
\nextref\Sk
\nextref\Vo

\topmatter

\abstract     We associate a contractible  ``outer space'' to any
free product of
groups $G=G_1*\dots*G_q$. It equals Culler-Vogtmann space
when $G$ is   free, McCullough-Miller space when no $G_i $ is $\Z$. Our
proof of contractibility (given   when $G$ is not free) is based on Skora's
idea of deforming morphisms between trees.

Using the action of $\Out(G)$ on this space, we show that
$\Out(G)$
      has finite virtual cohomological dimension, or is VFL (it has a finite
index subgroup with a finite classifying space), if the groups $G_i$
and $\Out(G_i)$ have similar properties. We deduce that $\Out(G)$ is
VFL if $G$ is a torsion-free hyperbolic group, or a limit group (finitely
generated fully residually free group).

\endabstract

\title The outer space   of a free product
           \endtitle

\author  Vincent Guirardel, Gilbert Levitt
     \endauthor

\endtopmatter

\document

\head \sect Introduction  and statement of results \endhead

A famous theorem of Grushko implies that a finitely generated group $G$
has a decomposition as a free product of
     subgroups $G=G_1*\dots*G_p*F_k$, with $F_k$   free
     of rank $k$ and every $G_i$ non-trivial,
non isomorphic to $\Z$, and freely indecomposable.
If $G=H_1*\dots*H_q*F_\ell$ is another such decomposition, the number
of factors is the same, $\ell=k$,  and (after reordering) $H_i$ is conjugate
to
$G_i$.

Despite this uniqueness, there is a lot of freedom in the choice of the
free factors, even when $k=0$ (but $p\ge3$). Because
of this, the  automorphism group
$\Aut(G)$ of $G$ is much more complicated than the direct product of
the groups $\Aut(G_i)$ and $\Aut(F_k)$. In particular, $\Aut(G)$
contains automorphisms acting on $G_i$ as conjugation by an
    element  of
$G_j$, with $j\neq i$, and acting on the other factors as the identity. A
presentation of
$\Aut(G)$ was given by Fouxe-Rabinovitch [\FR], generalizing work of
Nielsen on
$\Aut(F_n)$.

A different approach was introduced by Culler-Vogtmann [\CV], who
obtained finiteness results for $\Out(F_n)$ by letting it act on a
contractible complex now known as outer space (see [\Vo]).
McCullough-Miller [\McM] constructed a complex to study the group $\Sigma
\Aut(G)$ of symmetric automorphisms of a free product. This group equals
$\Aut(G)$ only when no factor is infinite cyclic, so in a sense these
two works cover opposite situations.

In this paper we construct a complex $P\O$ that allows studying
$\Aut(G)$ and $\Out(G)$ in all cases. We prove for instance:

    \thm{Theorem 5.2} Let $G=G_1*\dots*G_p*F_k$ be a Grushko
decomposition of a finitely generated group $G$. Assume that $Out(G)$ is
virtually torsion-free.
If $G_i$ and $\Aut(G_i)$ have finite virtual cohomological dimension
(resp\. have a finite index subgroup with a finite classifying space) then
so does
$\Out(G)$.
\fthm

Since a torsion-free finitely generated group is a free product of
cyclic groups and one-ended groups, knowledge of $\Out(H)$ for $H$
one-ended gives information about automorphisms of arbitrary
torsion-free groups.

We pointed out earlier that, in spite of the uniqueness in Grushko's
theorem, there is a lot of flexibility in free products. On the other hand,
one-ended groups often exhibit a strong algebraic rigidity. For instance,
Bowditch showed [\Bo] how to obtain the JSJ splitting of a one-ended
(word) hyperbolic group $H$ from the topology of its boundary. This
implies that the splitting is completely unique, not just up to certain
moves, and therefore invariant under automorphisms. Using this
invariance, one may describe $\Out(H)$ in terms of abelian groups and
mapping class groups (see [\Le] for a precise statement of this result of
Sela's [\Se]).

We deduce:

\thm{Theorem 6.1} If $G$ is a torsion-free hyperbolic group, then
$\Out(G)$ has a finite index subgroup with a finite classifying space. If
$G$ is a   virtually torsion-free hyperbolic group, then
$\Out(G)$ has  finite virtual cohomological dimension.
\fthm

Similarly,  a limit group (finitely generated
fully residually free group) has an invariant JSJ splitting if it
is one-ended, and one gets:

\thm{Theorem 6.5} If $G$ is a limit group, 
then
$\Out(G)$ has a finite index subgroup with a finite classifying space.
\fthm

Let us now explain how the complex $P\O$ is constructed. Consider a
Grushko decomposition
$G=G_1*\dots*G_p*F_k$ as above.
Like Culler-Vogtmann's outer space, $P\O$ may be viewed as a space
of $G$-trees, or as a space of marked   metric graphs. We describe it
here as a space of graphs. 

\midinsert
\centerline{\includegraphics[scale=1]{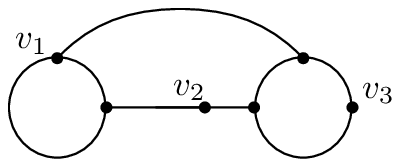}}
\botcaption{Figure 1}{A   graph of groups in
the outer space of
$G_1*G_2*G_3*F_3$}
\endcaption
\endinsert

A point in $P\O$ is defined by a marked  graph  of groups
$\Gamma$ as in fig 1. More precisely, $\Gamma$ is a finite graph of groups
with trivial edge groups; its edges are assigned a positive length;
the marking is an isomorphism from $G$ to $\pi_1(\Gamma)$, well defined
up to composition with inner automorphisms; for each $i\in\{1,\dots,p\}$
there is a vertex $v_i$ whose group is conjugate to $G_i$; all other vertex
groups are trivial and every terminal vertex is a $v_i$.

A point of $P\O$ is such a marked metric graph, projectivized (all edge
lengths may be multiplied by the same amount). The group $\Out(G)$
acts on $P\O$ by change of marking.  The quotient consists of  finitely
many open cells.

As always,  the difficulty is to show contractibility of the space. Our
approach is not combinatorial as in [\CV] or [\McM], but geometric. It
is based on the idea of deforming  morphisms
between trees, introduced   in Skora's unpublished paper [\Sk] giving a
geometric proof that Culler-Vogtmann space is contractible. We shall give
the construction   in its natural context, that of
$\bbR$-trees, but in this paper we only apply it to metric
simplicial trees.

Let   $f:T_0\ra T_\infty$ be a map between
metric simplicial trees, sending edges isometrically onto edges. We
factor it through intermediate trees $T_t$ ($t\ge0$) defined as follows:
two points $x,y\in T_0$ with $f(x)=f(y)$ are identified in $T_t$ if and only
if the image of the segment $[x,y]$ by $f$ is contained in the $t$-ball
around $f(x)$. The map $f$ factors through maps $f_t:T_t\to T_\infty$. (We
are grateful   to L\. Mosher for making us realize that this definition is
equivalent to Skora's.)

One may also visualize the trees $T_t$ as follows.
If $f$ is not an
embedding, then it folds two edges having a common vertex. Let $l$
be the length of the shortest pair of edges of $T_0$ folded by $f$.
Then, for $t\leq l$, the  tree $T_t$ is defined by folding along
length $t$ any pair
of edges that are folded by $f$.  For $t>l$ with   $t-l$ small enough, one
can similarly define
$T_t$ from $T_l$,   using the map $f_l$  and folding  along
     length $t-l$. It may   be shown that one can reach any value
$t\in\bbR^+$ by iterating this process finitely many times.

    This   deformation is
   constructed and studied in Section 3 (after the present paper was first
posted, we became aware of [\Cl], which contains another   account of
Skora's idea; see [\LL] for yet another  deformation).
Its most   delicate feature   is   continuity.

\thm{Theorem} The assignment $(f,t)\mapsto f_t$ is a continuous
semi-flow on the space of morphisms between \Rt s (equipped with
  the equivariant Gromov-Hausdorff topology).
\fthm

    In Section 4,
we define   our outer space and we   show:

\thm{Theorem} For any finitely generated group
$G=G_1*\dots*G_p*F_k$, with $p\ge1$, the outer space
$P\O$ is   contractible.
\fthm

This is proved as follows. We find   a basepoint $T_0\in P\O$
such that, for all
$T\in P\O$, one can define a morphism $f_T:T_0\ra T$ varying
continuously with $T$.  Applying the semiflow to $f_T$  (backwards) then
gives a continuous way of deforming $T$ into  $T_0$   (to be precise, we
work in the non-projectivized space $\O$, and
the space is contracted to a simplex rather than to a point, as the
edge lengths on $T_0$ depend on $T$).

As pointed out in [\McM], there are two
natural topologies on $P\O$. First, there is the weak topology associated
to its natural simplicial structure. Unlike Culler-Vogtmann space, the
complex $P\O$ is   usually not locally finite, and the weak topology is
different from the equivariant Gromov-Hausdorff topology, or axes
topology. The deformation argument gives contractibility in the axes
topology, whereas applications require weak contractibility. This is
obtained by showing that,   as $T_0$ and $T$ each vary within a given
simplex, the set of intermediate trees $T_t$   only meets finitely many
simplices.

The   main difference between our arguments and those used by
Skora in his proof that Culler-Vogtmann space is contractible [\Sk] is that
constructing the maps
$f_T$ is easier in our situation, because selecting the unique point of $T$
fixed by $G_1$ provides a continuous choice of a basepoint in $T$ as $T$
varies in $ \O$.

Because of this, contractibility of Culler-Vogtmann space is not proved
here. In [\GLdef], we    
prove contractibility of arbitrary deformation spaces, using a general
basepoint argument (two simplicial $G$-trees are in the same
deformation space if they have the same elliptic subgroups, see [\Fo]).
This applies in particular to    Culler-Vogtmann space, and to spaces of
JSJ splittings.

In Section 5, we obtain general finiteness results about $\Out(G)$ by
studying its action on $P\O$. In Section 6, this is applied to hyperbolic
groups, limit groups, and groups acting freely   on $\R^n$-trees.

\head \sect Trees \endhead

Trees will be considered both as combinatorial objects, and as metric
objects.

     Combinatorially, a (non-metric)
    simplicial tree  is a simplicial $1$-complex with no circuit. It is a
$G$-tree if the group $G$ acts on it by automorphisms, without
inversions. A map between simplicial trees is {\it simplicial\/} if
     it maps every edge onto an edge (in particular, no edge is collapsed).

An \Rt{} may be defined as a connected metric space   whose distance $d$
satisfies the $0$-hyperbolicity inequality
$$d(x,y)+d(z,t)\le\max(d(x,z)+d(y,t),d(x,t)+d(y,z))$$
(see [\AB]). Actions on \Rt s will always be by   isometries.

Let $T$ be a simplicial $G$-tree with finitely many $G$-orbits of edges.
It  becomes an \Rt{}
when each edge is assigned a positive length. An \Rt{} obtained in this
way will be called a {\it metric simplicial tree\/}. A metric simplicial
tree is thus determined by the underlying simplicial tree, and one
positive number for each orbit of edges.

The intermediate trees $T_t$ will be defined as \Rt s, but readers
uncomfortable with
\Rt s may assume all trees to be simplicial. In this
paper,  the semi-flow will only be used with $T_0$ and
$T_\infty$   simplicial. In this case, all intermediate trees $T_t$ are
simplicial.

The segment between two points $a$ and $b$ in a tree will be denoted
by $[a,b]$.
In an \Rt, the length of a segment $J=[a,b]$ will be denoted by $|J|$.
    A {\it finite subtree\/}  of a tree is the convex hull of a
finite set of points.

Let $T$ be a $G$-tree. The tree, or the action, is called
   {\it non-trivial\/} if there is no global fixed point,
{\it minimal\/}  if there is no proper invariant
subtree. In a minimal simplicial tree, there are only finitely many orbits
of edges if $G$ is finitely generated.

Maps between $G$-trees will always be assumed to be $G$-equivariant. A
map $f:T\to T'$  between \Rt s is a {\it morphism\/}  if
each segment in $T$ can be written as a  finite union of
subsegments, each of which is mapped
isometrically into $T'$ by $f$. If $T$ and $T'$  are metric
simplicial trees, $f$ is a
morphism if and only if one may subdivide
$T$ and $T'$ so that $f$  becomes simplicial and maps every edge isometrically.

We always assume  morphisms $f:T\to T'$ to be surjective. This is
automatically true when $T'$ is a minimal $G$-tree.

\head \sect The semi-flow \endhead

\subhead \subsect Definition of  the semi-flow \endsubhead

Let $f :T_0\to T_\infty$ be a  morphism between \Rt s. Our
goal is to construct for each $t\in \bbR^+$ an \Rt{} $T_t $,   and
morphisms
$\varphi _t :T_0\to T_t $ and $\psi _t :T_t \to T_\infty$, such that the
following diagram commutes:
$$\xymatrix@=2ex{ T_0 \ar[rr]^{f} \ar[rd]_{\phi_t }&&T_\infty\\ &T_t
\ar[ur]_{\psi_t } }$$
If $f $ is an equivariant morphism between two $G$-trees, then
$T_t$ is a
$G$-tree  and $\phi _t,\psi _t $ are equivariant.

When we need to study how this construction depends on $f$ and $t$,
we will write $T_t(f)$ instead of $T_t$, and $\Phi  (f,t),\Psi  (f,t)$
instead of $\varphi _t,\psi _t$.

\subsubhead The tree $T_t$ as a set \endsubsubhead

Let $d_0$ (resp\. $d_\infty$) denote   distance in $T_0$ (resp\.
$T_\infty$). Given $a,b\in T_0$ with $f (a)=f(b)$, we define their
identification time by
$$\tau (a,b)=\max_{x\in[a,b]} d_\infty\bigl(f(x),f(a)\bigr).$$     Note that
$\ds\tau (a,b)\le\frac{d_0(a,b)}2$.

For $t\ge0$, say that $a\sim_t b$ when $f(a)=f(b)$ and $\tau (a,b)\le
t$.   Obviously
$\tau (a,c)\le\max(\tau (a,b),\tau (b,c))$ if
$f(a)=f(b)=f(c)$, so that $\sim_t$ is an equivalence relation. We define
$T_t$   as $T_0/{\sim_t}$ and $\varphi _t :T_0\ra T_t$ as the quotient
map. The
      map $f$ induces a map
$\psi _t :T_t\to T_\infty$. The definitions may be extended to
$t=\infty$, by setting
$\varphi _\infty=f$ and $\psi _\infty=id_{T_\infty}$.

\subsubhead The metric on $T_t$ \endsubsubhead

This metric  will be the  maximal metric making $\varphi _t$
$1$-Lipschitz.

Let $a,b$ be  arbitrary points of $T_0$.  Given   $t\geq 0$,   a   {\it
$t$-subdivision\/} $\sigma $ between $a$ and $b$ is a sequence $
(a=y_{0},x_1,y_1,\dots,x_n,y_n, x_{n+1}=b)$ such that
$y_i\sim_t x_{i+1}$.
      When there is no risk of confusion, we will simply say
\emph{subdivision}. A subdivision  is {\it straight\/} if all points
belong to
$[a,b]$ and lie in the indicated order.

The {\it flesh\/} of $\sigma $ consists of the segments
$I_i=[x_i,y_i]$.  The \emph{jumps} of $\sigma$ are the segments
$J_i=[y_i,x_{i+1}]$. The {\it length\/} $|\sigma|$ of $\sigma $ is the
total length of its flesh:
$|\sigma |=\sum_{i=1}^n|I_i|$.

      Note that the concatenation of the paths $f(I_i)$  gives a path
joining $f(a)$ to $f(b)$ in $T_\infty$. This   implies the inequality
$|\sigma |\ge d_\infty(f(a),f(b))$. Similarly, the concatenation of the
paths $\phi_t(I_i)$ joins $\phi_t(a)$ to $\phi_t(b)$,    so to make
$\phi_t$ $1$-Lipschitz we need to choose a metric on $T_t$ so that  $d
(\phi_t(a),\phi_t(b))\le |\sigma|$. The metric $d_t$ which we take on
$T_t$ will be defined as
$d_{t}(\phi_t(a),\phi_t(b))=\inf_{\sigma}|\sigma|$, where the infimum
is taken over all $t$-subdivisions between $a$ and
$b$.  We first show that this infimum is achieved.

\nom\stra
\thm{Lemma \sta\ (Admissible subdivisions)} Given
$a,b\in T_0$ and
$t\in\bbR^+$,  there exists a straight $t$-subdivision $\sigma _0$
between
$a$ and $b$ such that $|\sigma _0|\le|\sigma |$ for every
$t$-subdivision
$\sigma $ between $a$ and $b$.
\fthm

Such a subdivision will be  called {\it admissible\/} (or $t$-admissible
if there is a risk of confusion).

\demo{Proof} All subdivisions considered  in this proof
are  $t$-subdivisions between $a$ and $b$. Define the  {\it
complexity\/} of a subdivision as
$c(\sigma )=(|\sigma |,
\sum_{i=0}^n|J _i|)$, with  $\sum |J_i|$  the  total length of the  jumps.
Complexities are compared in lexicographic order.

Let $K$ be  a finite subtree of $ T_0$ containing $a$ and $b$. We say
that a subdivision is contained in $K$ if all points in the sequence are
in $K$. We first show that    there exists a number $M_K$ such that,
given
$\sigma
\inc K$, there exists $\sigma '\inc K$ consisting of at most $M_K$
points, with
$c(\sigma ')\le c(\sigma )$.

Since $f$ is a morphism, we may cut $K$
into finitely many arcs, all of which are mapped injectively into
$T_\infty$. If two non-consecutive points of $\sigma $ belong to the
same arc, we may decrease $n$
         by removing points in-between. This does not   increase
        complexity, because    $|\mu  |\ge d_\infty(f(x),f(y))$ for any
subdivision $\mu $ between $x$ and $y$, and shows the existence of
$M_K$.

Since $\sim_t$ is a closed equivalence relation, a simple compactness
argument on the set of subdivisions of cardinality   at most $M_K$
shows that, among all subdivisions contained in
$K$, we can find $\sigma _K$ with minimal complexity.

We now  prove that $\sigma _K$ is straight.  Suppose it is not.
Exchanging $a$ and $b$ if needed, we may assume that there is an
overlap between $[x_i, y_i]$ and $[y_i,x_{i+1}]$. Since $x_{i+1}$ and
$y_i$ have the same image $p$ by $f$, and
$f$ is a morphism, the set $(y_i,x_{i+1}]\cap f\mi(p)$ is nonempty. Let
$z$  be its point  closest to $y_i$. We have $z\sim_t x_{i+1}$, because
$\tau (z,x_{i+1})\le\tau (y_i,x_{i+1})\le t$.

Note that the image of
$(y_i,z)$ by
$f
$ is contained in one component of $T_\infty\setminus\{p\}$.
Denoting by $y_{i,\varepsilon}$ (resp\. $z_\varepsilon $)    the point of
$(y_i,z)$ at distance
$\varepsilon $ from $y_i$ (resp\. $z$), this
implies
$y_{i,\varepsilon} \sim_t z_\varepsilon $ for $\varepsilon >0 $ small
enough. For such an $\varepsilon $,
consider
$\sigma _\varepsilon
$ obtained from $\sigma _K$  by changing   $(\dots x_i, y_i,
x_{i+1}\dots)$ into $(\dots x_i, y_{i,\varepsilon},z_\varepsilon ,z,
x_{i+1}\dots)$. We have $|\sigma _\varepsilon |=|\sigma _K|$, but
$c(\sigma _\varepsilon )<c(\sigma _K)$, a contradiction. This shows
that $\sigma _K$ is straight.

It is now easy to conclude. Define $\sigma _0$ minimizing
complexity among all subdivisions contained in $[a,b]$. Given any
subdivision $\sigma $, choose $K$ containing $\sigma $. Then
$|\sigma |\ge|\sigma _K|$, and $|\sigma _K|\ge|\sigma _0|$ because
$\sigma _K$ is straight.
\cqfd\enddemo

        We  define $\delta _t(a,b)$ as $|\sigma  |$, with $\sigma $
      a $t$-admissible subdivision between $a$ and $b$. Clearly,
$\delta_t(a,b)=0$ if and only if $a\sim_t b$. Moreover, $\delta_t$ is a
pseudo-metric on $T_0$. It follows that $\delta_t$ defines a metric
$d_t$ on $T_t$, with $d_t(\varphi _t(a),\varphi _t(b))=\delta _t(a,b)$.

Note that $d _0$ is just the original distance on $T_0$. Given $a$ and
$b$, the function $t\mapsto
\delta _t(a,b)$ is non-increasing and equals $\delta
_\infty(a,b)=d_\infty(f(a), f(b))$  for $\ds t\ge \frac{d_0(a,b)}2$.

Consider an admissible subdivision $\sigma$ between $a$ and $b$. Its
flesh   arcs map isometrically into $T_t$, and the concatenation of
their images is a path of
$d_t$-length
$|\sigma|=\delta _t(a,b)$  joining     $\phi_t(a)$ to $\phi_t(b)$. This
path is therefore geodesic. The fact that
$T_t$ is an
\Rt{} will imply that this is \emph{the} arc joining $\phi_t(a)$ to
$\phi_t(b)$.

\thm{Lemma \sta} The metric space $T_t$ is an \Rt.
\fthm

\demo{Proof}  Since $T_t$ is connected, it suffices to show that the
distances between any four points satisfy the $0$-hyperbolicity
inequality [\AB, Theorem 3.17]. If not, we can find $\theta >0$, and
$a_0,b_0,c_0,d_0\in T_0$ satisfying the following equation
$I_\theta $:
$$\delta _t(a_0,b_0)+\delta _t(c_0,d_0)=\theta +\max\bigl(\delta
_t(a_0,c_0)+\delta _t(b_0,d_0),
\delta _t(a_0,d_0)+\delta _t(b_0,c_0)\bigr) .$$

Fix such   $\theta, a_0,b_0,c_0,d_0$. Among all quadruples
$Q=(a,b,c,d)$ contained in the convex hull of $\{a_0,b_0,c_0,d_0\}$ and
satisfying $I_\theta $, choose one for which the total length of the
convex hull $C(Q)$ is minimal. We may assume that $a$ is a terminal
point of $C(Q)$.  Consider three admissible subdivisions
$(a, x_1,y_1,\dots)$ between
$a$ and the other three points.

Each of these subdivisions has positive length, because $I_\theta $
forces
        $a,b,c,d$ to have   distinct images in $T_t$. We may therefore
assume that the first flesh interval $I_1=[x_1,y_1]$   is
non-degenerate in all three subdivisions.  Then $x_1\neq a$ in at least
one subdivision,  since otherwise we could decrease the length of
$C(Q)$ by moving $a$ (without losing the equation
$I_\theta $). Thus $a\sim_t a' $ for some $a'\neq a$ in
$C(Q)$. Since $a',b,c,d$ satisfies $I_\theta$, one can decrease the total
length of
$C(Q)$, a contradiction.
\cqfd\enddemo

\thm{Lemma \sta}  The maps $\varphi _t:T_0\to T_t$ and $\psi
_t:T_t\to T_\infty$ are   morphisms.
\fthm

\demo{Proof} They are obviously surjective, since $f$ is surjective
(being a morphism). Now let
$I$ be any arc   on which
$f$ is isometric. Since
$\varphi _t$ and
$\psi _t$ are both $1$-Lipschitz,
       $\phi_t$ is isometric in restriction to   $I$, and  $\psi _t$ is
isometric in restriction to   $\varphi _t(I)$. This implies that
$\phi_t$ is a morphism. Furthermore $\psi _t$ is also a morphism,
since any arc of
$T_t$ is contained in a finite union of images  $\varphi _t(I)$.
\cqfd\enddemo

Assume that $T_0$ and $T_\infty$ are endowed with isometric actions
of a group $G$, and   $f$ is equivariant. It is clear that $T_t$ inherits an
isometric action of $G$ and that
$\phi_t$ and $\psi_t$ are equivariant. We will see  that $T_t$ is
simplicial if $T_0$ and
$T_\infty$ are, but minimality of
$T_0$ and $T_\infty$ as $G$-trees does not imply minimality of
$T_t$ (see 3.3).

\subhead \subsect Continuity \endsubhead

Fix a  discrete group $G$.  Let $\A$ be the space of $G$-trees, i.e\.
\Rt s $T$ with an isometric action of $G$. Distance will always be
denoted by $d$. Two trees are considered equal if there is an equivariant
isometry between them. The set $\A$ is equipped with the {\it
equivariant    Gromov-Hausdorff topology\/}.

Recall that a fundamental system of neighborhoods for $T\in\A$ is
given by the sets
$V_T(X, A,\varepsilon )$, with $X\inc T$  and
$A\inc G$    finite sets, and  $\varepsilon >0$. By definition,  $T'$ is in $
V_T(X, A,\varepsilon )$ if and only if there exists a ``lifting'' map
$x\mapsto x'$ from
$X$ to $T'$ such that $|d(x,gy)-d(x',gy')|<\varepsilon $ for every
$x,y\in X$ and $g\in A$.

Now consider the space $\M$ of (surjective)  $G$-equivariant
morphisms
$f:T\to T_\infty$ between
$G$-trees. Its topology is defined by neighborhoods $W_f(X,
A,\varepsilon )$,  with $f':T' \to T'_\infty$ in
$W_f(X, A,\varepsilon )$ if and only if there exists $x\mapsto x'$ as
above, with the extra requirement
$|d(f(x),f(gy))-d(f'(x'),f'(gy'))|<\varepsilon $. Note  that we define the
same topology if we drop $g$ from the last requirement, and that the
source and target maps from
$\M$ to
$\A$ are  continuous.

The construction given   in the previous section associates
morphisms $\varphi _t=\Phi (f,t)$ and $\psi _t=\Psi (f,t)$ to any
$f\in\M$ and $t\in[0,\infty]$.

\nom\continu
\thm{Proposition \sta} The maps $\Phi ,\Psi :\M\times [0,\infty]\to\M$
are continuous.\fthm

The main step in the proof is the following lemma:

\nom\contin
\thm{Lemma \sta}  Fix $f\in\M$, $\varepsilon >0$,
$t\in[0,\infty)$, and  $a,b\in T$. There exist $\alpha >0$ and a finite set
$F\inc [a,b]$ containing $a$ and $b$,   such that, if
$f'\in W_f(F,\{1\},\alpha )$ and
$|s-t|<\alpha $, then $|\delta _{s}(a',b')-\delta _t(a,b)|<\varepsilon $.
\fthm

In this statement, $a'$ and $b'$ are lifts of $a,b$ to the source
       of $f'$ provided by the definition of $W_f(F,\{1\},\alpha )$, and
$\delta _{s}$ is computed with respect to
$f'$.

The lemma immediately implies the continuity of $\Phi $ on
$\M\times  [0,\infty)$.  Continuity for $t=\infty$ is clear because, if
$f'\in W_f(\{a,b\} ,\{1\},\alpha )$, then $\delta _{s}(a',b')$ is constant
for $t\ge
\frac12(d(a,b)+\alpha )$. Continuity of $\Psi  $ is also clear once we
know that of $\Phi $, so \continu{} follows  from \contin. There remains
to prove Lemma \contin.

      \demo{Proof of Lemma \contin} In the following proof,
$C,C_1,C_2 $ will denote universal constants (which could be made
explicit).

      Since $\delta _t(a,b)$ is determined by the restriction of $f$ to
$[a,b]$, we may assume
$T=[a,b]$. We may also forget about $G$, so we simply write
$W_f(F, \alpha )$.  Let $E$ consist of
$a$, $b$, and all points where
$f$ folds. Consider $f'\in W_f(E ,\alpha )$.

Let $x\mapsto \ov x$ be the linear map taking
$[a,b]$ to $[a',b']$ (so $\ov a=a'$ and $\ov b=b'$).  It distorts distances
by at most
$\alpha$ (by which we mean $|d(x,y)-d(\ov x,\ov y)|\leq \alpha$). We
claim that
$$|d(f'(\ov x),f'(\ov y))-d(f(x),f(y))|\leq C \alpha \ntag$$ for all
$x,y\in[a,b]$.
      Such an inequality clearly holds for points of $E$ because  $\ov x$ is
$(C_1\alpha)$-close to $x'$ for $x\in E$.
      It is also true for arbitrary points because
$x\mapsto f'(\ov x)$ distorts distances by at most $(C_2\alpha )$ on
any interval on which $f$ does not fold. We shall assume $C\ge1$.

We first prove the lemma when $a\sim_t b$ and no point of $(a,b)$ has
the same image in $T_\infty$ as $a$ and $b$. In this case we take
$F=E$. Choose
$\theta$ such that the $\theta $-ball  around $f(a)$ contains no vertex
of the finite tree $f([a,b])$ different from $f(a)$, and no image of a
point where
$f$ folds.  Choose $\alpha $   with $10C\alpha <\varepsilon $ and
$4C\alpha <\theta $. Consider $f'\in W_f(F ,\alpha )$ and $s\ge t-\alpha
$.

Viewing $[a,b]$ as a subinterval of $\R$, we  have $f(a+ \beta
)=f(b-\beta )$ for $0\le \beta
\le \theta  $.  The intervals $I=[a+2C\alpha ,a+4C\alpha ]$ and
$J=[b-4C\alpha ,b-2C\alpha ]$ are mapped isometrically onto the
same arc in $T_\infty$.  Using $(1)$, we can then find
$a_1\in\rond{I}$ and $b_1\in \rond{J}$ with
$f'(\ov a_1 )=f'(\ov b_1)$. Note that $f([a_1,b_1])$ is disjoint from
$f([a,a+2C\alpha ])$. In particular, since $\tau (a,b)\le t$, every point
of $f([a_1,b_1])$ is
$(t-2C\alpha )$-close to $f(a_1)$. By $(1)$, we have $\tau (\ov a_1 ,\ov
b_1)\le t-2C\alpha+C\alpha  \le s$ and therefore
$\ov a_1
\sim_{s}\ov b_1$. We deduce $$\delta _{s}(\ov a,\ov b)\le \delta
_{s}(\ov a,\ov a_1)+
\delta _{s}(\ov a_1,\ov b_1)+
\delta _{s}(\ov b_1,\ov b)\le 10C\alpha <\varepsilon $$ since $\delta
_{s}(\ov a,\ov a_1)\le d(\ov a,\ov a_1)\le d(a,a_1)+\alpha \le 4C\alpha
+\alpha \le5C\alpha $. The lemma is proved in this case.

When $a\sim_t b$ but $A=f\mi(f(a))$ is larger than $\{a,b\}$, we take
$F=E\cup A$ and we prove the lemma by applying the previous
argument on each subinterval of $[a,b]$  bounded by points of
$A$.

In the general case, we fix a  $t$-admissible subdivision
$\sigma =(a,x_1,y_1,\dots,x_n,y_n,b)$ between $a$ and $b$. Define
$F$ as the union of   of $E$ and all preimages $f\mi(f(x))$, for $x$ a
point of
$\sigma
$.  Let $f'\in W_f(F ,\alpha )$ and $s\in(t-\alpha ,t+\alpha )$. We will
force $|\delta _{s}(\ov a,\ov b)- \delta _{t}(a,b)|<\varepsilon $ by
choosing
       $\alpha $ small enough.

We have $\delta _t(a,b)=\sum_{i=1}^nd(x_i,y_i)$, and
$$\delta _{s}(\ov a,\ov b)\le \sum_{i=1}^nd(\ov x_i,\ov y_i)
+\sum_{i=0}^n\delta _{s}(\ov y_i,\ov x_{i+1}).$$ Furthermore, $d(\ov
x_i,\ov y_i)\le d(x_i,y_i)+\alpha $. Since $y_i\sim_t x_{i+1}$,
      we have seen that we can choose $\alpha$ small enough so that
$\ds\delta _{s}(\ov y_i,\ov x_{i+1})\le \varepsilon /{2(n+1)}$. If
furthermore $n\alpha <\varepsilon /2$ we obtain
$\delta _{s}(\ov a,\ov b)\le \delta _{t}(a,b)+\varepsilon $.

We now get a lower bound for $\delta _{s}(\ov a,\ov b)$. We may
assume that $f$ does not fold on $I_i=[x_i,y_i]$ (if $t=0$, we may have
to refine $\sigma $).  The segment $[(\ov x_i)_s,(\ov y_i)_s]\inc T'_s$
then has length at least
$|I_i|-C\alpha $. As above, we choose $\alpha $ so that
$(\ov y_i)_s$ and
$(\ov x_{i+1})_s$ are $(\varepsilon /2(n+1))$-close.

      We shall prove that the intersection between $[(\ov x_i)_s,(\ov
y_i)_s]$ and $[(\ov x_{i+1})_s,(\ov y_{i+1})_s]$ has length at most
$3C\alpha $. For $\varepsilon $ and $C\alpha $ small with respect to
$\ds\min_i|I_i|$ we can then write
$$ \delta _{s}(\ov a,\ov b)\ge\sum _{i=1}^n(|I_i|-C\alpha
)-\sum_{i=0}^n(\frac\varepsilon {2(n+1)}+3C\alpha )\ge \delta
_t(a,b)-4(n+1)C\alpha -\frac\varepsilon 2.
$$ This gives the required inequality provided $\ds\alpha
<\frac\varepsilon {8(n+1)C}$.

If the intersection is bigger than $3C\alpha $, we can find $\ov u\in
[\ov x_i,\ov y_i]$ and $\ov v\in [\ov x_{i+1},\ov y_{i+1}]$ with $\ov
u_s=\ov v_s$, and both $d(\ov u, \ov y_i)$,
$d(\ov v, \ov x_{i+1})$ bigger than $3C\alpha $. We shall reach a
contradiction by considering the corresponding points $u\in I_i$ and
$v\in I_{i+1}$. Let
$p=f(y_i)=f(x_{i+1})\in T_\infty$.

We have $d( f(u),p)=d(u,y_i)\ge3C\alpha -\alpha \ge2C\alpha $.
Similarly,
$d(f(v),p)\ge2C\alpha $. On the other hand,  $d(f(u),f(v))\le C\alpha
$.  This implies that the arcs $f(I_i)$ and $f(I_{i+1})$  belong to (the
closure of) the same component $\C$ of $T_\infty\setminus \{p\}$. In
particular, the points
$y_i-\kappa $ and $x_{i+1}+\kappa $ have the same image in
$T_\infty$ for
$\kappa >0$ small. Since they don't have the same image in $T_t$,
there exists
$z\in[y_i,x_{i+1}]$ such that $f(z)$ is at distance exactly $t$ from
$p$, and in a component of
$T_\infty\setminus \{p\}$ other than  $\C$. This implies that the
distance from $f(z)$ to $f(u)$ is at least $t+2C\alpha $, and therefore
by $(1)$ the point $f'(\ov z)$ has distance at least $t+C\alpha $ from
$f'(\ov u)=f'(\ov v)$. But
$\ov u_s=\ov v_s$, so $s\ge t+C\alpha \ge t+\alpha $, a contradiction.

This completes the proof of Lemma \contin, hence also of Proposition
\continu. \cqfd
\enddemo

\subhead \subsect Additional properties of the semi-flow
\endsubhead

The material from this section will not be used in an essential way in the
rest of the paper, as alternative arguments will be provided (see
Remark 4.1).

Let $f:T_0\to T_\infty$ be a  fixed (surjective) morphism between
$G$-trees, with $G$ finitely generated.

\subsubhead Simplicial   trees
\endsubsubhead

\nom\simpli
\thm{Proposition \sta} If
$f:T_0\to T_\infty$ is a morphism   between metric simplicial
$G$-trees
      with finitely many orbits of edges, then the tree $T_t$ is simplicial
for every $t>0$.
\fthm

\demo{Proof}  After subdividing, we may assume that each edge of
$T_0$ is mapped isometrically onto an edge of $T_\infty$. Fix $t$. Let
$S_\infty\inc T_\infty $ consist of all vertices, and all points lying at
distance exactly $t$ from some vertex. Let
$S_0$ and $S_t $ denote the preimages of $S_\infty$ in $T_0$ and
$T_t$. Since there are finitely many orbits of edges in $T_\infty$, all
three sets $S_0$,  $S_t $,
$S_\infty$ are closed and meet a given finite subtree in finitely many
points.

We now show  that any point $u\in T_t$ not in $S_t$ is regular: its
complement in $T_t$ has only two components. This implies that $T_t$ is
simplicial, with vertex set contained in $S_t$.

Consider $x,y\in T_0$ mapping onto $u$. They belong to open edges of
$T_0$. The identification time $\tau (x,y)$ is at most $t$. If $\tau (x,y)=t$,
consider $p$ at distance $t$ from
$f(x)$  in $f([x,y])$.  It is a vertex of $T_\infty$, since $f$ is an
immersion away from vertices. This implies  $\psi _t(u)\in S_\infty$,
contradicting $u\notin S_t$. We deduce $\tau (x,y)<t$, and there exist
open neighborhoods of $x$ and $y$ in $T_0$ with the same image in
$T_t$.  It follows that $u$ is  regular.
\cqfd\enddemo

\subsubhead Minimality \endsubsubhead

Let $T$ be an \Rt{} with an isometric action of
$G$. Recall that   $T$ (or the action) is called minimal if there is no proper
$G$-invariant
     subtree. If the action of $G$ on $T$   is   non-trivial (there is no global
fixed point), there is a unique minimal  subtree
$\Tmin\inc T$; it consists of all  axes of  hyperbolic elements of $G$
(see [\CuMo]).

A {\it terminal point\/} of   $T$ is a point $x$ which is not contained in
any  open interval (equivalently, $T\setminus\{x\}$ is   connected). If
$T$ is minimal, it contains no terminal point (because
$T\setminus G.x$ is a proper invariant subtree if $x$ is terminal).

The action of $G$ on    $T$ is {\it finitely supported\/} if there exists
a finite subtree   $K$ such that every segment  of $T$ is covered by
finitely many translates of $K$.  A simplicial action is finitely
supported if and only if there are  finitely many   orbits of edges.
A minimal action  of a  finitely
     generated group   on an \Rt{} is   finitely supported.

\thm{Lemma \sta}
A non-trivial finitely
supported $G$-tree $T$   is minimal if and only if it has no terminal
point.
\fthm

\demo{Proof} We have seen that there is no terminal point in a minimal
tree. We now assume that $T$ is not minimal and we find a
terminal point.
Let $\Tmin$ be the  minimal subtree.
Consider  a finite tree $K$ such that   $G .K=T$. Then $K\cap
\Tmin$ is a convex subset of
$K$, distinct from $K$, so there is at least one endpoint of $K$ which is
not contained in $\Tmin$. Consider an endpoint $x$ of $K$ whose
distance $d$ to $\Tmin$ is maximal. Clearly, any point of $T$ is at
distance at most $d$ from $\Tmin$. This implies that $x$ is terminal in
$T$.
\cqfd\enddemo

    Let $f:T_0\to T_\infty$ be as above.
Minimality of
$T_0$ and $T_\infty$ as $G$-trees does not imply minimality of
$T_t$. For instance, if $T_0$, $T_\infty$ and $f$ are simplicial, and
$f$ maps all   edges incident to some vertex $v\in T_0$ onto the same
edge of $T_\infty$,   then  $\phi_t(v)$ is a terminal vertex of $T_t$ for
$t$ small enough.
This is essentially the only way   in which  $T_t$   may  fail
to be minimal.

    \thm{Definition \sta} A morphism of $\bbR$-trees $f:T_0\ra T_\infty$
{\rm  satisfies the minimality condition} if,  for all
$x\in T_0$, there exists an open interval containing $x$ on which $f$ is
one-to-one. In particular, $T_0$ and
$T_\infty$
      have no terminal point, and are therefore minimal if  they are
finitely supported.
\fthm

\example{Remark \sta} One can easily prove   that, if $f:T_0\ra
T_\infty$ is a  morphism between  minimal simplicial
$G$-trees, then
      one can construct another morphism $g:T_0\ra T_\infty$ satisfying
the minimality condition.
\endexample

The minimality condition clearly implies minimality of $T_t$:

\nom\mini
\thm{Proposition  \sta} Let
$f:T_0\to T_\infty$ be a morphism   between $G$-trees. Assume that
$T_0$ is finitely supported  and $f$ satisfies the minimality
condition.  Then for all $t\ge0$ the tree $T_t$ is finitely  supported, and
the morphisms $\phi_t
$, $\psi_t $ satisfy the minimality condition. In particular,   $T_t$ is
minimal. \cqfd
\fthm

\subsubhead The semi-flow property \endsubsubhead

\thm{Proposition  \sta} For all $s,t\in\bbR$, we have $\Psi((\Psi
(f,t),s)=\Psi (f, s+t)$ and  $\Phi(\Psi(f,t),s)\circ\Phi(f,t)=\Phi(f,s+t)$:
folding $T_t$ along distance $s$  gives the same result as folding
$T_0$ along
$s+t$.
\fthm

$$\xymatrix@=2ex{T_0\ar[rrr]^f \ar[dr]_{\Phi(f,t)} &&&T_\infty \\
            &T_t \ar[dr]_{\Phi(\Psi(f,t),s)} \ar[urr]|{\Psi(f,t)}\\
            && T_{s+t} \ar[uur]_{\Psi(\Psi(f,t),s)=\Psi(f,{s+t})}  }
$$

\demo{Proof} Fix
$s,t\ge0$, and denote by $a_t$ the image  of a point $a\in T_0$ in
$T_t=T_t(f)$. We have to show $\delta _{s+t}(a,b)=\delta _s(a_t,b_t)$
(with $\delta _s$ computed with respect to $\Psi (f,t):T_t\to
T_\infty$).
      We are going to prove the equality
$$\tau (x_t,y_t)  =\max(0,\tau (x,y)-t), \ntag$$  valid whenever
$f(x)=f(y)$ (with $\tau (x_t,y_t) $ also   computed with respect to
$\Psi  (f,t)$).

We first check that (2) implies the lemma. By (2), the projection onto
$T_t$ of an $(s+t)$-subdivision $\sigma$ between $a$ and
$b$ is an $s$-subdivision $\sigma_t$  between $a_t$ and $b_t$. Since
$|\sigma_t|\leq |\sigma|$, one gets  $\delta _{s+t}(a,b)\geq\delta
_s(a_t,b_t)$. Conversely,  consider an $s$-subdivision
$\sigma_t$ between
$ a_t$ and $b_t$ in $T_t$, and lift it arbitrarily to   $\sigma$ between
$a$ and $b$. By (2), $\sigma $ is an $(s+t)$-subdivision. If $[x_i,y_i]$
are the arcs in the flesh    of
$\sigma $, then $\delta _{s+t}(a,b)\le\sum_i\delta _{s+t}(x_i,y_i)\le
\sum_i \delta _{ t}(x_i,y_i)
      = |\sigma _t|$. This implies $\delta _{s+t}(a,b)\le\delta _s(a_t,b_t)$.

Let us now prove (2).  Fix a $t$-admissible subdivision $\sigma_0 $
between $x$ and $y$.  The set
$A=\Psi  (f,t)([x_t,y_t])\inc T_\infty$  is the image by
$f$ of the flesh of
$\sigma_0 $. To get the image  of $[x,y]$ by $f$, one has to add the
image of the jumps $[y_i,x_{i+1}]$. These are loops based at points of
$A $, each contained in the $t$-ball around its basepoint. This shows
$\tau (x,y)\le\tau (x_t,y_t)+t$, because $A$ is contained in the $\tau
(x_t,y_t)$-ball around $f(x)$.

We conclude  by proving the opposite inequality $\tau
(x_t,y_t)\le\tau (x,y)-t$, when $\tau (x,y)\ge t$ (note that
$x_t=y_t$ when $\tau (x,y)\le t$). Let $B_\infty\subset T_\infty$ be
the set of points
$p\in  f([x,y])$ such that the segment $[f(x),p]$ may be extended
      geodesically by a segment $[p,q]\inc f([x,y])$ of length $t$. Note that
$B_\infty$ is contained in the
$(\tau (x,y)-t)$-ball around $f(x)$, so we need only  prove
$\Psi  (f,t)([x_t,y_t])\subset B_\infty$.

So consider $B_t$, the preimage of $B_\infty$ in $[x_t,y_t]$, and
assume that
$B_t\subsetneqq [x_t,y_t]$.  The set $B_t$ is closed,   and contains
$\{x_t,y_t\}$ because $\tau (x,y)\ge t$. Consider a
      component $(u_t,v_t)$ of $[x_t,y_t]\setminus B_t$, and choose
preimages $u_0,v_0 $    of
$u_t,v_t$ in $[x,y]\inc T_0$, with
$[u_0,v_0]$  minimal  (for   inclusion) among all possible choices. Note
that $f(u_0)=f(v_0)$ and $\Phi (f,t)((u_0,v_0))\subset (u_t,v_t)$.
Since $f((u_0,v_0)) =\Psi  (f,t)((u_t,v_t)) $
     does not meet   the component of  $T_\infty\setminus\{f(u_0)\}$
that contains $f(x)$, it is contained in the $t$-ball around
$f(u_0)$. We get $u_0\sim_t v_0$ and $u_t=v_t$, a contradiction.
\cqfd\enddemo

\head \sect The outer space  \endhead

We fix a finitely generated group $G$,
decomposed as $G=G_1*\dots*G_p*F_k$ with each
$G_i$ non-trivial. We require $p\ge 1$ and $p+k\ge2$, but in  this
section $G_i$ may be $\Z$ or a nontrivial free product. In particular, the
space constructed in [\McM] to study symmetric automorphisms is a
special case of the space $P\O$ defined below.

All trees  considered here will be minimal simplicial
$G$-trees, up to equivariant isometry. All maps between trees will be
$G$-equivariant (hence onto). Unless otherwise indicated, we
     assume that there is no  {\it redundant vertex\/}:   if
$v$ has degree
      $2$, then $v$ is the fixed  point of an element of $G$
      exchanging the two edges incident to $v$.

Let $\O$ be the space of  {\it metric simplicial $G$-trees\/} $T$ such
that:

$\bullet$ The action of $G$ on $T$ is minimal, with trivial edge
stabilizers.

$\bullet$ For each $i$, there is exactly one orbit of vertices with
stabilizer conjugate to $G_i$.

$\bullet$ All other points have trivial stabilizer.

Via Bass-Serre theory, an element of $\O$ may also be viewed as a
{\it marked metric graph of groups\/}. It is a finite graph of groups
$\Gamma
$, with an isomorphism from
$\pi _1(\Gamma )$ to
$G$ (the marking), well-defined up to composition with inner
automorphisms. Edge groups are trivial. There is one vertex
$v_i$ with group  conjugate    to
$G_i$ for each
$i$,   and all other vertex groups are trivial. Edges are assigned a
positive length.   Minimality of the action on $T$ translates to the fact
that every terminal vertex is a $v_i$.

The quotient of $\O$ by the natural action of $(0,\infty)$ (defined by
rescaling, i.e\. multiplying all lengths by the same number)  is the {\it
projectivized space\/}
$P\O$ (the outer space of $G$).
We equip $\O$ with the equivariant  Gromov-Hausdorff topology, and
$P\O$ with the quotient topology.

The action of $G$ on $T$ defines a length function $\ell_T:G\to\R$, by
$\ds\ell_T(g)=\min_{x\in T}d(x,gx)$. As  actions in $\O$ are minimal
with non-abelian length function, assigning to $T\in\O$ its length
function defines  homeomorphisms from $\O$ onto its image in
$\R^G$, and from $P\O$ onto its image in $P\R^G$ [\Pa].

We often
identify $\O$ and $P\O$ with their images. The topology on $\O$ and $P\O$ is
   commonly called
    the length function topology, or  the {\it axes topology\/}.

We now consider the simplicial structure on $\O$ and $P\O$.  The
metric on  a tree $T\in\O$ is determined  by finitely many positive
numbers, one length for each orbit of edges (equivalently, one length
for each edge of $\Gamma $). The set of trees obtained from
$T$ by varying these numbers will be called the (open simplicial) {\it
cone\/} containing $T$. In other words, a cone consists of equivariantly
   homeomorphic trees (trees with the same underlying
simplicial tree).

      The projection  of a  cone to $P\O$ is an {\it open simplex\/}. Its
closure in $P\R^G$ is a closed simplex $\Sigma $, whose points may be
viewed as projectivized length assignments   on   the set of edges of
$\Gamma $. The intersection of $\Sigma
$ with $P\O$ will be called a {\it closed simplex\/} of $P\O$ (a point of
$\Sigma $ is in
$P\O$ if and only if
      the union of   edges of length $0$ is a forest, and every component
of the forest   contains at most one
$v_i$, so   a closed simplex of
$P\O$ is a  finite union of open simplices).

Just like Culler-Vogtmann's outer space, $P\O$ is not a simplicial
complex: simplices may have faces ``at infinity''. If one insists on
having a simplicial complex, one replaces $P\O$ by its barycentric
spine (there is an equivariant deformation retraction of  $P\O$ onto its
spine, as  in [\CV,
\McM]).

The {\it weak topology\/} on $P\O$ is defined in the usual way: a set is
closed if and only if its intersection with every closed simplex is
closed. This topology does not coincide with the axes
topology, because the simplicial structure is not locally finite (see
[\McM]). The two topologies have the same restriction, however, on
any finite union of simplices (see [\GLdef] for a detailed study).

Our goal now is to show that $P\O$ is contractible, both in the axes
topology and in the weak topology. Our main tool will be the semi-flow.

\nom\stable
\example {Remark \sta} $\O$ is invariant under the semi-flow, in
the following sense. Let
$f:T_0\to T_\infty$ be a morphism between trees in $\O$, and $t>0$.
The tree $T_t$ is simplicial by  Proposition \simpli{} (this will also
follow from Lemma 4.3).  It clearly has the correct edge and vertex
stabilizers, because both $T_0 $ and $T_\infty$ do.

Furthermore, $T_t
$ is minimal if $f$ maps each edge isometrically and all vertex stabilizers
of
$T_0$ are non-trivial. This follows from Proposition \mini, but here is a
direct argument. Given an edge $vw$ of $T_0$, choose a non-trivial
element
$g_v$ (resp\. $g_w$) fixing
$v$ (resp\. $w$). The   translation axis of $g_vg_w$ in $T_t$ contains the
images of
$v$ and $w$, so the minimal subtree of $T_t$ contains every vertex and
$T_t$ is minimal.
\endexample

\nom\contra
\thm{Theorem \sta} The   space $P\O$  is contractible in the axes
topology.
\fthm

\demo{Proof} We use the semi-flow to construct a contraction in
$\O$. We first define a basepoint (or rather a basecone) in $\O$ (see
figure 2).

\midinsert
\centerline{\includegraphics[scale=1]{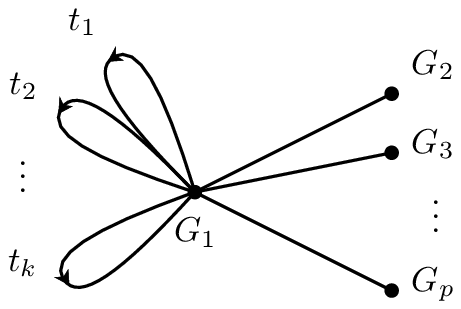}}
\botcaption{Figure 2}{The base point of outer space}
\endcaption
\endinsert

Choose a free basis
$t_1,\dots, t_k$ of $F_k$, and let $T_0\in \O$   be described by Figure 2.
Formally, $T_0$ is such that every vertex has non-trivial stabilizer, and the
vertex stabilized by $G_1$ is adjacent to the vertices stabilized by
      $G_2,\dots, G_p, t_1G_1t_1\mi,\dots, t_pG_1t_p\mi $.
    This uniquely defines
$T_0$ as a simplicial   non-metric $G$-tree, as we haven't specified
edge-lengths
yet.

Given $T \in\O$, we define an equivariant map
$f:T_0\to T $ by mapping every vertex of $T_0$ to the unique vertex
of $T $ with the same stabilizer, and extending linearly on every edge.

There is a unique edge-length assignment on
$T_0$ such that
$f$ is isometric on every edge (in particular, $f$ is a morphism). We
thus associate to $T\in\O$ a metric tree $T_0(T)\in\O$ and a morphism
$f_T:T_0(T)\to T$. As a simplicial tree, $T_0(T)$ does not depend on
$T$   and thus lies in the cone consisting of all metric
simplicial trees obtained
by varying lengths on $T_0$.

Let $T_t(T)$ be the tree defined by applying the semi-flow to
$f_T$. It belongs to $\O$ (see Remark \stable).

We can now consider the map $\rho :(T,t)\mapsto T_t(T)$, from
$\O\times[0,\infty] $ to $\O$.   Once we know that $\rho$ is
continuous,  we will get a continuous family of paths  connecting $T$
to
$T_0(T) $, thus providing a   deformation   retraction of $\O$ onto a
cone. This will prove   contractibility of $\O$.

Let us check  that $\rho $ is continuous (in the axes
topology).  The   edge-length assignment on $T_0$ depends
continuously on $T$: if $g_i\in G_i$ is non-trivial, the edge of
$\Gamma $ between $v_1$ and $v_i$ has length
$\frac12\ell_T(g_1g_i)$, and the loop associated to $t_j$ has length
$\frac12\ell_T(g_1t_jg_1t_j\mi)$.  More generally, if $x,y$ are
vertices of $T_0$ and $g_x,g_y$ are non-trivial elements fixing them,
then $d(f_T(x),f_T(y))=\frac12\ell_T(g_xg_y)$. It follows that
$T_0(T)$ and $f_T $ depend continuously on $T$, so $\rho $ is
continuous because  the semi-flow is continuous
(Proposition
\continu).

To prove contractibility of $P\O$,  we simply reparametrize
$\rho
$ so that it descends to $\hat\rho :P\O\times[0,\infty] \to P\O$. For
instance, we may take  $\wtilde\rho (T,t)= \rho (T,\ell_T(h)t)$, with
$h\in G$ chosen so that
$\ell_T(h)>0$ for every $T\in\O$ (take $h=g_1g_2$ if $p\ge2$ and
$h=g_1t_jg_1t_j\mi$ if
$p=1$, with 
$g_i \in G_i$   non-trivial).
\cqfd\enddemo

To prove contractibility in the weak topology, we simply show that
$\hat\rho$ is also continuous in the weak topology. This  requires the
following  ``finiteness lemma''.

\nom\finit
\thm{Lemma \sta}   Consider two   non-metric simplicial
$G$-trees $T_0$, $T$,
with finitely many orbits of edges. Let $f:T_0\to T$ be an
equivariant simplicial map (sending   every edge onto an edge)
such that the preimage
$f\m(e)$ of each  edge $e\inc T$
      is a finite set of edges.

If $f$ factors through  equivariant continuous maps
$$\xymatrix@=2ex{T_0\ar^f[rr]\ar[dr]_{f_1}&&T\\&T'\ar[ur]_{f_2}}$$
    with $T'$ an \Rt{} and  $f_1$  surjective,
then $T'$ is a simplicial tree, and  there are only finitely many
possibilities for   $T'$    up to  equivariant homeomorphism (i.e\. as a
non-metric simplicial
$G$-tree).
\fthm

\example{Remark}
In this lemma, $T_0$ and $T$ are allowed to have redundant vertices.
\endexample

\demo{Proof}  For each open edge $e$ of $T$, let  $K_{e,i}$ ($i\in I_e$)
be the (finite) set  of connected components of the closure of
$f_2\mi(e)$.
Each $K_{e,i}$ is a
finite subtree of $T'$. If $(e,i)\neq(e',j)$, then  either
$K_{e,i}\cap K_{e',j}=\ev$, or $K_{e,i}\cap K_{e',j}$ consists of one
point, which is the image of an endpoint of an edge in $f\m(e)$. In
particular,  for a given $K_{e,i}$, the set of points occurring as
$K_{e,i}\cap K_{e',j}$ is finite.

It follows that $T'$ is simplicial.  Furthermore,  if a vertex $v$ of $T'$ is
not the image of a vertex of $T_0$, then there exist   
edges $e_0,e'_0$  of
$T_0$, with the same image in $T$, such that $v$ is an endpoint of
$f_1(e_0)\cap f_1(e'_0)$. This  implies that there is a uniform bound
(depending only on $f$) for the number of $G$-orbits of vertices of
$T'$.

Now   subdivide the trees so as to make $f_1$ and
$f_2$ simplicial: we first subdivide $T $ by adding images of vertices
of
$T'$, then we  subdivide $T_0$ and $T'$ by adding all preimages of
vertices of $T $.  This may create redundant  vertices, but the number
of orbits of redundant vertices in
$T$ and
$T_0$ is uniformly bounded; in particular,  there are only finitely many
possibilities for the subdivisions of $T$ and $T_0$.

We fix these subdivisions, and  we show that there are only finitely
many possibilities for $T'$. Let   
    $e_1, \dots, e_n$ be  a set of 
representatives for $G$-orbits of oriented edges of $T$ (subdivided). In
the set $E$ of  oriented edges of (the subdivided) $T_0$, define $e\sim e'$ 
if $e,e'$ have the same image in $T'$. There are only finitely many
possibilities for $\sim$, since it is determined by its restriction to the
finite set $\cup_i f\mi(e_i)$.

We conclude the proof by showing
that $\sim$ completely determines $T'$   (as a non-metric
$G$-tree). Consider the connected graph $T_0/\sim$ obtained from $T_0$ by
identifying
$e,e'$ whenever $e\sim e'$.
Since $T'$ is obtained from
$T_0/\sim$ by identifying vertices, but no edges,  simple connectedness of
$T'$  implies that $T'=T_0/\sim$.
\cqfd\enddemo

\nom\contrweak
\thm{Corollary \sta} The space $P\O$  is contractible in the weak
topology.
\fthm

\demo{Proof}
We first deduce from Lemma \finit{} that, if $S$ is an open
simplex of
$P\O$, then the restriction of $\hat\rho $ to
    $S\times[0,\infty]$ meets only finitely many simplices.
     Let $\Tilde S$ be the preimage of $S$ in $O$.
Given $T\in \Tilde S$, consider the map $f_T$
constructed in the proof of Proposition \contra{}  and subdivide
$T_0$
    so that
$f_T$ becomes simplicial (this creates redundant vertices). As
a simplicial map (forgetting edge lengths),
$f_T$ does not depend on   $T\in \Tilde S$.
    Two distinct edges of $T_0 $ with
the same image in $T $ are in different
$G$-orbits (because $T   $ has trivial edge stabilizers), so
the
finiteness condition in Lemma \finit{} is satisfied. The lemma then
implies that for  $T\in\Tilde S$ there are only finitely many
possibilities for $T_t(T)$
as a simplicial tree.

We can now prove that $\hat\rho $ is continuous in the weak topology.
Since a closed simplex $\Sigma $ of $P\O$ is the union of finitely many
open simplices, the restriction of $\hat\rho $ to
$\Sigma \times[0,\infty]$ is continuous in the weak  topology because it
is continuous  in the axes topology, and the two topologies
coincide on any finite union of simplices.  By definition of the weak
topology, this implies continuity of $\hat\rho $.
\cqfd\enddemo

\remark{Remark}   In \cite{\GLdef}, we   extend these results and
prove the  contractibility of   any \emph{deformation
space}. In particular, this  includes outer space and    spaces of JSJ
splittings of   finitely presented groups. This extended result
also implies that the  fixed point set   of any finitely generated  
subgroup
$F$ of
$\Out(G)$ acting  on $P\O$ is contractible or empty (and non-empty if
$F$ is finite and   solvable).
\endremark

\head \sect The action of $\Out(G)$\endhead

We denote by $Z(H)$ the center of a group $H$. The group of inner
automorphisms is $Inn(H)\simeq H/Z(H)$, and the group of outer
automorphisms is $\Out(H)=\Aut(H)/Inn(H)$. Note that a nontrivial free
product has trivial center.

Let $G=G_1*\dots*G_p*F_k$ be as above, with $p\ge 1$ and $p+k\ge2$.
We now assume that $G_i$ is freely indecomposable and not isomorphic
    to
$\Z$. In particular, $G$ is a nontrivial free product and is not free.

Since any  automorphism of $G$ maps $G_i$ onto a conjugate of some
$G_j$, there is a natural action of
$\Out(G)$ on the contractible complex $P\O$, obtained by
precomposing   actions on trees by automorphisms. If an element of
$P\O$ is viewed as a graph of groups $\Gamma
$, the action is by changing the marking.

Let $S$ be an open simplex of $P\O$. We will think of it as a simplicial
$G$-tree $T$ (with no metric specified), or as a marked graph of
groups $\Gamma
$. If we forget the marking of $\Gamma $, there are only finitely many
possibilities for $\Gamma $ (because of the minimality assumption).
This means that  there are  only finitely many orbits of simplices
under the action of $\Out(G)$ on $P\O$.

Let us now describe the stabilizer of $S$, that is the subgroup
$\Out^S(G)\inc\Out(G)$ sending
$S$ to itself. It is the group of   automorphisms preserving the
decomposition of
$G$ as a graph of groups given by $\Gamma $. Such groups have been
studied in [\BJ] and [\Le].

Let $\Out_0^S(G)\inc \Out^S(G)$ be the finite index subgroup
consisting of automorphisms acting trivially on the quotient graph
$\Gamma =T/G$. Since  the edge groups of $\Gamma $ are trivial,
$\Out_0^S(G)$ has a very simple description (see [\Le, Prop.\ 4.2]).
Denoting   by $n_i$   the degree of the vertex $v_i$ in $\Gamma $, it 
 is a direct product
$ \prod_{i=1}^p M_{n_i}(G_i)$ of groups which fit in exact sequences
$$\{1\}\to G_i^{n_i}/Z(G_i)
\to M_{n_i}(G_i)\to\Out(G_i)\to\{1\}\ntag$$
       $$\{1\}\to G_i^{n_i-1}
\to M_{n_i}(G_i)\to\Aut(G_i)\to\{1\},\ntag$$   with the center
     $Z(G_i)$   embedded diagonally into
$G_i^{n_i}$ (the groups    $M_{n_i}(G_i)$ are denoted by
$\MCG^{\partial} (G_i)$ in  [\Le], as they are ``pure mapping class groups'').

The exact sequence $(4)$   is split and $M_{n_i}(G_i)$ is
the semi-direct product $G_i^{n_i-1}
\rtimes\Aut(G_i)$ associated to the diagonal action of $\Aut(G_i)$ on
$G_i^{n_i-1} $. In particular, $M_{1}(G_i)=\Aut (G_i)$, and
$M_{2}(G_i)$ is the holomorph Hol$(G_i)$.

Let $\Out'(G)\inc\Out(G)$  be the finite index subgroup
consisting of automorphisms mapping each
$G_i$ to a conjugate of itself. It maps onto
$\Out(G_i)$ in  a natural way (since $G_i$ equals its normalizer in
$G$).  We also consider the homomorphism $\Out'(G)\to\Out(F_k)$
obtained by viewing $F_k$ as the quotient of $G$ by the normal
subgroup generated by the
$G_i$'s.

\nom\triang
\thm{Lemma \sta} If all groups $G_i$ and $G_i/Z(G_i)$ are
torsion-free, then the kernel of $\pi :
\Out'(G)\to\Out(F_k)\times\prod_i\Out(G_i)$ is torsion-free.
\fthm

\demo{Proof}  If not, choose $\Phi \in\ker\pi $ whose order is a prime
number
$p$. It fixes a point in $P\O$, since otherwise $\Z/p\Z$ would act
freely on
$P\O$ and have a finite-dimensional classifying space. This means that
$\Phi $ belongs to some $\Out^S(G)$.

Consider the  action of $\Phi $ on the quotient graph of groups
$\Gamma $ corresponding to $S$.  It acts trivially on the topological
fundamental group $\pi _1(\Gamma )\simeq F_k$, and sends each
vertex $v_i$ to itself.   Thus $\Phi $ acts trivially on $\Gamma $, so
$\Phi
\in\Out_0^S(G)$.

We have seen that $\Out_0^S(G)= \prod_{i }
M_{n_i}(G_i)$. Since $\Phi $ maps trivially into
$\Out(G_i)$, it belongs to the product $\prod_i G_i^{n_i}/Z(G_i)$. But
$H^{n }/Z(H)$ is  torsion-free if both $H$   and
$H/Z(H)$ are (if $h_1^p=\dots=h_n^p\in Z(H)$  with $p\ge2$, then $h_i\in
Z(H)$ because
$H/Z(H)$ is torsion-free, and $h_i=h_j$ because
$(h_ih_j\mi)^p=h_i^ph_j ^{-p}=1$).  
\cqfd\enddemo

We can now prove:  

\nom\applic
\thm{Theorem \sta }
\roster
\item  "{(i)}" Assume that $G_i$ and $G_i/Z(G_i)$ are torsion free, and
    $\Out(G_i)$ is virtually torsion-free. Then
$\Out(G )$ and $\Aut(G)$ are  virtually torsion-free.

If, furthermore,   $G_i$, $G_i/Z(G_i)$, and  $\Out(G_i)$ have finite  virtual
cohomological dimension (resp\. have a finite index subgroup  with a
finite classifying space), then so do  $\Out(G)$ and $\Aut(G)$.

\item "{(ii)}"   If  $G_i$ and $\Aut(G_i)$ have finite  virtual
cohomological dimension (resp\. have a finite index subgroup  with a
finite classifying space), then so do  $\Out(G)$ and $\Aut(G)$,
provided that
they are virtually torsion free.

\endroster
\fthm

    Basic facts about virtual cohomological
dimension and finite classifying spaces are recalled in   section 6 of
[\McM].

\remark{Remark}
    If $p+k\ge3$, then $G_i$, $G_i/Z(G_i)$ and $\Aut(G_i)$ are isomorphic to
subgroups of
$\Out(G)$. In particular, if $\Out(G)$ is virtually torsion
free (resp\. has finite cohomological dimension), so are $G_i$,
$G_i/Z(G_i)$ and $\Aut(G_i)$.
\endremark

\demo{Proof} Since $\Aut(G)$ is an
extension of $\Out(G)$ by
$\text{Inn}(G)\simeq G$,   it suffices to prove the required results for
$\Out(G)$.

     The first assertion of (i)   follows from Lemma
\triang, since $\Out(F_k)$ is virtually torsion-free.

To prove  the second assertion, we consider the action of
$\Out(G)$   on the contractible complex $P\O$ (strictly speaking, one
should replace $P\O$ by its barycentric spine as in [\McM]). We have
seen that there are only finitely many orbits of simplices, so it
suffices to check that all stabilizers $\Out ^S(G)$ have the property
under consideration.
As explained above,    $\Out ^S(G)$ is virtually a direct product of
groups
$M_{n_i}(G_i)$.  The exact sequences $(3)$ and $(4)$
then
imply that
$M_{n_i}(G_i)$, hence also $\Out ^S(G)$, has the required
property.

The proof of (ii) is similar. Note that the
     groups $M_{n_i}(G_i)$ are
       virtually torsion-free, because they are subgroups of
$\Out(G)$.
\cqfd\enddemo

Here is a slightly stronger result:

\nom\plusfort
\thm{Corollary \sta} Suppose that each $G_i$ is finitely generated and
has a normal subgroup of finite index
$H_i$, with
$H_i$ and $H_i/Z(H_i)$  torsion-free, and $\Out(H_i)$  virtually
torsion-free. Then $\Out(G)$ is virtually torsion-free. If furthermore
$H_i$ and $\Aut(H_i)$  have finite virtual cohomological dimension,
then
$\Out(G)$ has finite virtual cohomological dimension.
\fthm

The proof requires  the following fact (see [\Kr], [\McC]).

\nom\krs
\thm{Lemma \sta} Let $\Gamma $ be a finitely generated group. Let
$\Delta
$ be a normal subgroup of finite index,   with trivial center.
\roster
\item Some finite index subgroup $\Aut_0(\Gamma )\inc\Aut(\Gamma
)$ embeds into
$\Aut(\Delta )$.
\item  Some finite index subgroup $\Out_0(\Gamma
)\inc\Out(\Gamma )$ is isomorphic to the quotient of a subgroup of
$\Out(\Delta )$ by a finite subgroup.
\endroster
\fthm

\demo{Proof } Let $\Aut_0(\Gamma )$
        consist  of all automorphisms of $\Gamma$ mapping $\Delta $ to
$\Delta
$ and inducing the identity on $\Gamma /\Delta $, and let
$\Out_0(\Gamma )$ be the image of
$\lambda :\Aut_0(\Gamma )\to\Out(\Gamma )$.

We claim that the natural map from
$\Aut_0(\Gamma )$ to $\Aut(\Delta )$ is injective.   Assume that
$\alpha $ is in the kernel. Given $h\in \Gamma $, define the element
$\delta _0(h)\in\Delta $ by
$\alpha (h)=h \delta_0(h)$. For all $g \in\Delta$, we  have $h g h
\mi=\alpha (h g h \mi)=h \delta _0(h )g \delta _0(h)\mi h \mi$, so
$\delta _0(h )\in Z(\Delta )=\{1\}$ and $\alpha $ is the identity.

The group $\Sigma =\Aut_0(\Gamma )/Inn(\Delta )$ is therefore
isomorphic to a subgroup of
$\Out(\Delta )$. The kernel of $\lambda $ contains   $Inn(\Delta)$ with
finite index, so  $\Out_0(\Gamma )$ is the quotient of $\Sigma $ by a
finite subgroup.
\cqfd\enddemo

\nom\krss
\thm{Corollary \sta} Let $\Gamma $ be a finitely generated group. Let
$\Delta
$ be a normal subgroup of finite index,   with trivial center. If
$\Out(\Delta )$ is virtually torsion-free (resp\. has finite virtual
cohomological dimension), so   is $\Out(\Gamma )$. \cqfd
\fthm

\demo{Proof of Corollary
\plusfort}
       Let $\Delta $ the kernel of the map from $G$ to $\prod_i G_i/H_i$. It
is a finite free product of groups which are   isomorphic to $\Z$ or to
an $H_i$. Corollary  \plusfort{} now follows from Theorem
\applic {} and Corollary
\krss.
\cqfd\enddemo

\head \sect Applications\endhead

\subhead \subsect Hyperbolic groups\endsubhead

\thm{Theorem \sta} Let $G$ be a hyperbolic group.
\roster
\item If $G$ is torsion-free, then $\Out(G)$ has a finite index subgroup
with a finite classifying space.
\item If $G$ is virtually torsion-free, then $\Out(G)$ has finite virtual
cohomological dimension.
\endroster
\fthm

It is not known whether there exist hyperbolic groups which are not
virtually torsion-free (see [\KW]).

\demo{Proof} First assume that $G$ is torsion-free. There are two
cases (besides
$G=\{1\}$ and $G=\Z$).  If $G$ is one-ended, then some finite index
subgroup of $\Out(G)$ maps onto a direct product of mapping class
groups of  surfaces, with kernel $\Z^n$ (see [\Le], [\Se]).  The result is
therefore true in this case.  If
$G$ has infinitely many ends, we simply apply Theorem
\applic. Note that
$G_i$ has trivial center, and a finite classifying space (the quotient of
the Rips complex).

Assertion $(2)$ follows from Corollary \krss.
\cqfd\enddemo

\subhead \subsect   Groups acting freely on
$\R^n$-trees\endsubhead

In this section, we consider a   finitely generated group $G$
admitting a free action  on
an
$\R^n$-tree. This includes
limit groups  (also   known as finitely
generated
$\omega
$-residually free groups); this is due to   Remeslennikov [\Remes], see
[\Gu] for a correct statement.

The hypothesis that $G$
acts freely    on
an
$\R^n$-tree implies that $G$ is torsion-free and
commutative transitive (centralizers of nontrivial   elements are
abelian; see for instance [\Ch]). In particular, if
$G$ is not abelian, it has trivial center. Furthermore, $G$ has a finite
classifying space, every abelian subgroup is contained in   a maximal
one, maximal abelian subgroups  are finitely generated  and malnormal
(see [\Gr,
\Gu], and [\kmI,  \SeIHES] for limit groups).

We first   prove:

\nom\limi
      \thm{Theorem \sta} If a finitely generated group $G$ acts freely on
an
$\R^n$-tree, then $\Out(G)$ has finite virtual cohomological
dimension.
\fthm

If $n=1$, then $G$ is a free product of surface groups and free abelian
groups by Rips's theorem (see  [\BF], [\GLP]). Since Theorem
\limi{} is true for such groups, it is true for
$G$ by Theorem
\applic. The proof in the general case will be   by induction on $n$, using
the following fact.

\nom\recur
      \thm{Theorem \sta}  If a finitely generated, freely indecomposable
group $G$ acts freely on an
$\R^n$-tree ($n\ge2$), there exists an $\Out(G)$-invariant $G$-tree
$T$ with edge stabilizers  isomorphic to $\Z$, and vertex stabilizers
acting freely on
$\R^{n-1}$-trees. No edge   of $T$ joins two vertices with abelian
stabilizers.
\fthm

\demo{Proof}  
The tree $T$ is  constructed in Proposition 4.2 of [\Pab], as (a small
modification of) the cyclic JSJ splitting of
$G$.  The only property not proved there is that vertex stabilizers act
freely on $\R^{n-1}$-trees. But a vertex stabilizer $G_v$    is a
  surface group  (and therefore acts freely on an \Rt) or fixes a point in
every
$G$-tree
$T'$ with cyclic edge stabilizers   [\Pab, Th\'eor\`eme 4.1].    By
Theorem 7.1 of [\Gu],   there exists such a 
$T'$ with vertex stabilizers acting freely on
$\R^{n-1}$-trees. It follows that $G_v$ acts freely on an
$\R^{n-1}$-tree.
\cqfd\enddemo

Let $V$ be the vertex set of the graph of groups $\Gamma $
associated to $T$. By [\Le, Proposition 4.2], some finite index
subgroup $\Out_1(G)\inc\Out(G)$  fits in an exact sequence
$$1\to\T\to\Out  _1(G)
\to \prod_{v\in V } M(G_v)\to1,  
$$ where $G_v$ is the vertex  group, $M(G_v)$ is a subgroup of
$\Out(G_v)$, and
$\T$ is the {\it group of twists\/} associated to $\Gamma $.  The group
$\T$ is generated by a finite direct product  of centralizers of edge
groups in vertex groups (see [\Le], or the proof below). In particular,
$\T$ is a finitely generated abelian group.

The next proposition will show that $\T$ is torsion-free. Assuming
this, we complete the proof of Theorem \limi{} as follows.  

\demo{Proof of Theorem \limi }  
First
assume that $G$ is freely indecomposable. Apply Theorem \recur. By
the induction hypothesis, all groups $M(G_v)$ have finite virtual
cohomological dimension. Since
$\T$ is torsion-free with finite cohomological dimension,
$\Out_1(G)$ is virtually torsion-free and has finite virtual  cohomological
dimension. The same is
true of
$\Out(G)$. If $G$ is a non-trivial free product, we use Theorem
\applic.
\cqfd\enddemo

\thm{Proposition \sta} Let $\Gamma $ be a   minimal graph of groups
decomposition of a commutative   transitive group $G$, with edge
groups isomorphic to
$\Z$. Assume that every edge of $\Gamma $ has at least one endpoint with
nonabelian vertex group.
Then the group of twists $\T$ is isomorphic to a finite
direct product of abelian subgroups of $G$.
\fthm

\demo{Proof}  
 We denote by $E$ be the set of oriented edges of
$\Gamma $, by $E_v$     the set of edges with origin $v$, by $o(e)$  
the origin of $e $, by $Z_{G_{o(e)}}(G_e)$   the centralizer of the edge group
$G_e$ in the vertex group $G_{o(e)}$.  We say that a vertex
$v$ of $\Gamma $ is abelian if $G_v$ is  abelian. By subdividing, we   may
assume   that each edge joins an abelian vertex to a nonabelian one.

We recall the presentation of $\T$ given in [\Le,
Proposition 3.1].
The group
$\T$ is the quotient of $\ds\prod_{e\in E} Z_{G_{o(e)}}(G_e)$ by edge
and vertex relations, defined as follows.   For every
     pair
$(f,\ov f)$ of opposite edges, we kill  the diagonal image of
$G_f=G_{\ov f}$ in
$Z_{G_{o(f)}}(G_f)\times Z_{G_{o(\ov f)}}(G_{\ov f})\inc\prod
Z_{G_{o(e)}}(G_e)$.
     For every abelian vertex $v$, we kill the diagonal image of $G_v$ in
$\prod_{e\in E_v}  Z_{G_v}(G_e)\inc\prod _{e\in E} Z_{G_{o(e)}}(G_e)$  
(note that
$Z_{G_v}(G_e)=G_v$).
 There is no  relation at non-abelian vertices.

We first reduce to the case where every non-abelian vertex is
terminal, simply replacing a non-abelian
$v$ with  degree $n_v\ge2$  by $n_v$ terminal vertices each carrying
$G_v$. This may disconnect $\Gamma $, but then $\T$ is the direct
product of the groups associated to the components of the new graph.
We may therefore assume that   $\Gamma $ is a minimal graph of
groups  of the following form: it consists of a central abelian vertex $v$,
connected to
   terminal non-abelian vertices $v_1,\dots,v_n$ by edges
$e_1,\dots,e_n$.

First suppose that $G_{e_i}$ is equal to $Z_i$, its
centralizer in
$G_{v_i}$, for all $i$. Then $\T$ is the quotient of
$(G_v)^n$ by $G_v$ (embedded diagonally), so is isomorphic to
$(G_v)^{n-1}$. If $G_{e_i}$ is properly contained in $Z_i$  for   some $i$,
then by transitive commutativity $G_{e_i}$ maps onto
$G_v$. Furthermore, this may happen only for one value of $i$, say
$i=1$. The group $\T$ is the quotient of
$Z_1\times \prod_{i>1}  G_v$ by the image of $G_v$, embedded
diagonally into the whole product   (including the factor
$Z_1$). Since $G_{e_1}$ maps   onto
$G_v$, minimality of
$\Gamma $ implies $n\ge2$,    so $\T$ is isomorphic to $Z_1\times
\prod_{i>2}G_v$.
\cqfd\enddemo

\example{Remark}  The proposition is   true without the assumption on
edges of $\Gamma $, provided $G$ is not a solvable Baumslag-Solitar
group
$BS(1,s)$. 
It
also applies to abelian splittings,
provided $G$ is not the fundamental group of a graph of groups of the
following form: there is only one vertex $v$, its group
$G_v$ is abelian, and all inclusions from edge groups into $G_v$,
except possibly one, are onto.
\endexample

Our results so far used the  cyclic JSJ splitting of
$G$.  Using results of [\KM] and [\BKM] about the abelian JSJ splitting,
one gets:

\nom\vfl
\thm{Theorem \sta} If $G$ is a limit group, then $\Out(G)$ has a finite
index subgroup with a finite classifying space.
\fthm

\demo{Proof}  Since limit groups have finite classifying spaces, Theorem
\applic{} lets us assume that $G$ is one-ended. By Theorems 3.13 and
3.17 of [\BKM], it has an $\Out(G)$-invariant abelian JSJ splitting. The
associated graph of groups $\Gamma $ has three types of vertices:
elementary vertices (whose vertex group is a maximal abelian
subgroup), surface vertices, rigid vertices. Every edge connects an
elementary vertex to a non-elementary vertex.

Since this splitting is $\Out(G)$-invariant, we  may  use the exact
sequence
$$1\to\T\to\Out  _1(G)
\to \prod_{v\in V } M(G_v)\to1
$$ of [\Le] as in the proof of Theorem \limi.

Since   not   every
$\Out(G_e)$ is finite, proving that $\Out  _1(G)$ has finite index in
$\Out  (G)$ now requires the following  argument, extending
Proposition 2.3 of [\Le]. An edge $e$ with $\Out(G_e)$ infinite
connects an elementary vertex $v$ to a rigid vertex $w$. Since only
finitely many outer automorphisms of $G_w$ extend to
automorphisms of $G$ by [\KM,   Theorem 11.1], the group
$M(G_w)$ has finite index in the image of $\rho _w$.  Furthermore,
every bitwist around $e$ is a twist because the normalizer of $G_e$ in
an abelian $G_v$ equals its centralizer.

      As before, the group $\T$ is finitely generated and abelian. If $v$ is
rigid, then $M(G_v)$ is finite by [\KM]. If
$v$ is a surface vertex, $M(G_v)$ is a surface mapping class group.
If $G_v$ is abelian, then $M(G_v)$ is the subgroup of
$\Aut(G_v)=GL(n,\Z)$ consisting of automorphisms equal to the
identity on some finite collection of subgroups; it has a finite index
subgroup with a finite classifying space. Since we know that
$\Out(G)$ is virtually torsion-free, we conclude that some finite
index subgroup has a finite classifying space.
\cqfd\enddemo

Theorem \vfl{} is actually   valid  for a broader class of groups,
including groups acting freely on $\R^n$-trees.   This will appear
elsewhere, as it requires a general construction of invariant abelian
splittings.

\Refs\widestnumber\no{99}
\refno=0

\bref \by R. Alperin, H. Bass \paper Length functions  of group actions
on $\Lambda $-trees, {\rm in ``Combinatorial group theory and
topology (S.M. Gersten, J.R. Stallings, ed.)''} \jour Ann. Math. Studies
111
\publ Princeton Univ. Press \yr1987 \endref

\bref  \by H. Bass, R. Jiang  \paper Automorphism groups of tree
actions and of graphs of groups\jour J. Pure Appl. Algebra \vol 112
\yr1996\pages 109--155
\endref

\bref\by M. Bestvina, M. Feighn\paper Stable actions of groups on real
trees\jour Inv. Math.\vol 121\yr1995\pages287--321\endref

\bref \by B. Bowditch \paper Cut points and canonical splittings of
hyperbolic groups\jour Acta Math. \vol180\yr1998\pages145--186\endref

      \bref \by I. Bumagin, O. Kharlampovich, A. Miasnikov
\paper Isomorphism problem for finitely generated fully residually
free groups \jour preprint
\endref

\bref \by I. Chiswell \book Introduction to
$\Lambda$-trees\publ World Scientific Publishing Co.\yr2001
\endref

\bref\by M. Clay \paper Contractibility of deformation spaces of
$G$-trees
\jour Alg. \& Geom. Topology\vol5\yr2005\pages1481--1503\endref

\bref  \by M. Culler, J. Morgan \paper Group actions on \Rt s\jour
Proc. London Math. Soc.
\vol 55
\yr1987\pages571--604 \endref

\bref  \by M. Culler, K. Vogtmann \paper Moduli of graphs and
automorphisms of  free groups\jour Invent. Math. \vol 84
\yr1986\pages91--119 \endref

\bref   \by M. Forester\paper Deformation and rigidity of simplicial
group actions on trees\jour Geom.
\& Topol.\vol 6\yr2002\pages 219--267\endref

\bref \by D.I. Fouxe-Rabinovitch \paper On the automorphism
group of free products I
\jour Rec. Math. [Mat. Sbornik] N.S.\vol8\yr1940\pages 265--276
\moreref \paper II\jour Rec. Math. [Mat. Sbornik] N.S. \vol
9\yr1941\pages183--220
\endref

\bref \by D. Gaboriau, G. Levitt, F. Paulin\paper Pseudogroups of
isometries of $\R$ and Rips' theorem on free actions on \Rt s\jour
      Isr. Jour. Math.
      \vol 87\yr 1994\pages 403--428
\endref

\bref \by S. Gross \paper Group actions on $\Lambda$-trees \jour PhD
thesis, Hebrew University of Jerusalem\yr1998
\endref

\bref \by V. Guirardel \paper Limit groups and groups acting freely on
$\R^n$-trees
\jour Geometry and Topology \vol8\yr2004\pages1427--1470
\endref

\bref \by V. Guirardel, G. Levitt \paper Deformation spaces of
 trees\jour   arXiv:math.GR/0605545  \endref

\bref\by I. Kapovich, D. Wise \paper The equivalence of some residual
properties of word-hyperbolic groups \jour J. Algebra
\vol223\pages562--583\yr 2000\endref

\bref \by   O. Kharlampovich, A. Myasnikov \paper
Irreducible affine varieties over a free group. I.
Irreducibility of quadratic equations and Nullstellensatz
\jour J. Algebra
\vol200\yr1998\pages472--516
\moreref \paper II.
Systems in triangular quasi-quadratic form and description of residually
free groups
\jour J. Algebra
\vol200\yr1998\pages517--570\endref

      \bref \by   O. Kharlampovich, A. Myasnikov
\paper Effective JSJ   decompositions \jour Contemp. Math.   
(Groups, languages, algorithms)
 \yr 2005\vol378\pages87--212
\endref

\bref\by S. Krsti\'c \paper Finitely generated virtually free groups
have finitely presented automorphism group\jour Proc. London Math.
Soc.\vol64\yr1992\pages49--69
\endref

\bref\by G. Levitt \paper Automorphisms of   hyperbolic
groups and graphs of groups \jour Geom.
Dedic. \vol 114\yr2005
\pages 49--70 
\endref 

\bref\by J. Los, M. Lustig\paper The set of train track representatives
of an irreducible free group automorphism is contractible \jour
preprint
\endref

\bref\by J. McCool\paper The automorphism groups of finite
extensions of free groups\jour Bull. London Math. Soc.\yr
1988\vol20\pages 131--135\endref

      \bref   \by D. McCullough, A. Miller \paper Symmetric automorphisms
of free products \jour Mem. Amer. Math. Soc.
\vol122\yr1996\endref

\bref \by F. Paulin \paper The Gromov topology on $\R$-trees
\jour Topology and its App.\vol 32 \yr 1989\pages 197--221\endref

\bref \by F. Paulin \paper Sur la th\'eorie \'el\'ementaire des groupes
libres\jour S\'eminaire Bourbaki 922, Ast\'erisque   \vol 294\yr  
2004\pages 363--402\endref

\bref \by V. Remeslennikov \paper {$\exists$}-free groups as groups with a
length function
\jour Ukrainian Math. J. \vol44\yr1992\pages733--738\endref

\bref \by Z. Sela \paper Structure and rigidity in (Gromov) hyperbolic
groups and discrete groups in rank $1$ Lie groups II\jour GAFA
\vol7\yr1997\pages561--593\endref

\bref\by Z. Sela\paper Diophantine geometry over groups I:
Makanin-Razborov diagrams\jour Publ. Math.
IHES\pages31--105\yr2001\vol93
\endref

\bref \by R. Skora\paper  Deformations of length functions in
groups\jour
      Preprint
\yr1989\endref

\bref\by K. Vogtmann \paper Automorphisms of free groups and
outer space
\jour Geom. Dedic. \vol94\yr2002\pages 1--31
\endref

\endRefs

\address  V.G.: Laboratoire \'Emile Picard, umr cnrs 5580, Universit\'e
Paul Sabatier, 31062 Toulouse Cedex 4, France.\endaddress\email
guirardel\@picard.ups-tlse.fr{}{}{}{}{}\endemail

\address  G.L.: LMNO, umr cnrs 6139, BP 5186, Universit\'e de Caen,
14032 Caen Cedex, France.\endaddress\email
levitt\@math.unicaen.fr{}{}{}{}{}\endemail

\enddocument